\def\VersionDateTime{16/January/2014. Version 1.0.}
\documentclass[12pt, leqno]{amsart}
\usepackage{amscd, amsmath, amssymb, amsfonts}
\usepackage{latexsym}
\usepackage{mathrsfs}

\newtheorem{Theorem}{Theorem}[section]
\newtheorem{Proposition}[Theorem]{Proposition}
\newtheorem{Lemma}[Theorem]{Lemma}

\newtheorem{Claim}{Claim}[Theorem]

\theoremstyle{definition}
\newtheorem{Definition}[Theorem]{Definition}
\newtheorem{Remark}[Theorem]{Remark}

\renewcommand{\theTheorem}{\arabic{section}.\arabic{Theorem}}
\renewcommand{\theClaim}{\arabic{section}.\arabic{Theorem}.\arabic{Claim}}
\renewcommand{\theequation}{\arabic{section}.\arabic{equation}}

\newcommand{\KK}{{\mathbb{K}}}

\newcommand{\ZZ}{{\mathbb{Z}}}
\newcommand{\QQ}{{\mathbb{Q}}}

\newcommand{\OO}{{\mathcal{O}}}

\newcommand{\Xscr}{{\mathscr{X}}}
\newcommand{\Dscr}{{\mathscr{D}}}

\newcommand{\Div}{\operatorname{Div}}

\newcommand{\Gal}{\operatorname{Gal}}
\newcommand{\ord}{\operatorname{ord}}
\newcommand{\Pic}{\operatorname{Pic}}

\newcommand{\rest}[2]{\left.{#1}\right\vert_{{#2}}}  % restriction of #1 to #2
\renewcommand{\setminus}{\smallsetminus}

\newcommand{\Spec}{{\operatorname{Spec}}}

\newcommand{\outdeg}{\operatorname{outdeg}}

\newcommand{\Supp}{\operatorname{Supp}}
\newcommand{\alg}{\operatorname{alg}}

\newcommand{\Bs}{\operatorname{Bs}}

\newcommand{\Proof}{{\sl Proof.}\quad}
\newcommand{\QED}{{\unskip\nobreak\hfil\penalty50\quad\null\nobreak\hfil
{$\Box$}\parfillskip0pt\finalhyphendemerits0\par\medskip}}
\newcommand{\val}{\operatorname{val}}

\newcommand{\ch}{\operatorname{char}}

  {\end{list}}

  {\end{list}}
\begin{document}

%%%%%%%%%%%
%% Title              %%
%%%%%%%%%%%
\title{
Algebraic rank on hyperelliptic graphs and 
graphs of genus $3$}
\author{Shu Kawaguchi}
\address{Department of Mathematics, Graduate School of Science,
Kyoto University, Kyoto 606-8502, Japan}
\email{kawaguch@math.kyoto-u.ac.jp}
\author{Kazuhiko Yamaki}
\address{Institute for Liberal Arts and Sciences, 
Kyoto University, Kyoto, 606-8501, Japan}
\email{yamaki@math.kyoto-u.ac.jp}
\date{\VersionDateTime}
\thanks{The first named author partially supported by KAKENHI 24740015, and the second named author partially supported by KAKENHI 21740012.}

\newcommand{\Proj}{\operatorname{\operatorname{Proj}}}
\newcommand{\Prin}{\operatorname{Prin}}
\newcommand{\Rat}{\operatorname{Rat}}
\newcommand{\zero}{\operatorname{div}}

\begin{abstract}
Let $\bar{G} = (G, \omega)$ be a vertex-weighted graph, and $\delta$ a divisor class on $G$. Let $r_{\bar{G}}(\delta)$ denote the combinatorial rank of $\delta$. 
Caporaso has introduced the algebraic rank $r_{\bar{G}}^{\alg}(\delta)$ of $\delta$,  by using nodal curves with dual graph $\bar{G}$. 
In this paper, when $\bar{G}$ is hyperelliptic  or 
of genus $3$, we show that $r_{\bar{G}}^{\alg}(\delta) \geq r_{\bar{G}}(\delta)$ holds, generalizing our previous result. 
We also show that, with respect to the specialization map from a non-hyperelliptic curve of genus $3$ to its reduction graph, any divisor on the graph lifts to a divisor on the curve of the same rank. 
\end{abstract}

\maketitle

%%%%%%%%%
%  Introduction  %
%%%%%%%%%
\section{Introduction}
Let $k$ be an algebraically closed field. 
The correspondence between nodal curves over $k$ and their (vertex-weighted) dual graphs appears naturally in algebraic geometry, as in the description of the stratification of the Deligne--Mumford moduli space of stable curves. 
Recently, a theory of divisors on graphs has been developed (see for example \cite{AC}, \cite{BHN}, \cite{BN} and \cite{BN2}). This enables one to study the relationship between {\em linear systems} on a nodal curve and {\em those} on the corresponding graph (and also between linear systems on the generic fiber and those on the dual graph of the special fiber of a semi-stable curve over a discrete valuation ring): See, for example, \cite{ABBR}, \cite{B}, \cite{Ca}, \cite{Ca1}, \cite{Ca2}, \cite{CDPR} and \cite{KY}.  In particular, a tropical proof of the Brill--Noether 
theorem has been obtained in \cite{CDPR}. 

In this development, Caporaso \cite{Ca2} has defined the algebraic rank $r_{\bar{G}}^{\alg}(\delta)$ of a divisor class $\delta$ on a vertex-weighted graph $\bar{G} = (G, \omega)$, by using nodal curves with dual graph $\bar{G}$. 
It was shown in \cite[Summary~3.4]{Ca2} that, on some simple graphs $\bar{G}$, the algebraic rank $r_{\bar{G}}^{\alg}(\delta)$ equals the combinatorial rank $r_{\bar{G}}(\delta)$ for any divisor class $\delta$. 
Further, Caporaso, Len and Melo in \cite{CLM} have recently shown that 
$r_{\bar{G}}^{\alg}(\delta)  \leq r_{\bar{G}}(\delta)$ holds for any 
divisor class $\delta$ on any vertex-weighted graph $\bar{G}$. 
In  \cite[Proposition~1.5]{KY}, we showed that, if  $\ch(k) \neq 2$ and $\bar{G}$ is a hyperelliptic vertex-weighted graph satisfying a certain assumption 
on the bridges of $G$, 
then $r_{\bar{G}}^{\alg}(\delta)  \geq r_{\bar{G}}(\delta)$ holds for any 
divisor class $\delta$. 

In this paper, firstly, we show, based on \cite[Proposition~1.5]{KY},
that the above assumption on  
the bridges for hyperelliptic graphs is not necessary. 

\begin{Theorem}
\label{thm:cor:main}
Assume that $\ch(k) \neq 2$. 
Let $\bar{G} = (G, \omega)$ be a hyperelliptic vertex-weighted graph.  
Then, for any divisor class $\delta$ on $G$, we have 
$
  r_{\bar{G}}^{\alg}(\delta)  \geq r_{\bar{G}}(\delta). 
$
\end{Theorem}

Secondly, we show the same inequality on non-hyperelliptic graphs of genus $3$. 

\begin{Theorem}
\label{thm:main:2}
Let $\bar{G} = (G, \omega)$ be a vertex-weighted graph.  
Assume that $\bar{G}$ is non-hyperelliptic and of genus $3$. 
Then, for any divisor class $\delta$ on $G$, we have 
$
  r_{\bar{G}}^{\alg}(\delta)  \geq r_{\bar{G}}(\delta). 
$
\end{Theorem}

These results, combined with the above result of Caporaso--Len--Melo, show  that the algebraic rank equals the combinatorial rank on all hyperelliptic vertex-weighted graphs (when $\ch(k) \neq 2$) and non-hyperelliptic vertex-weighted graphs of genus $3$ (and certain graphs which are built from hyperelliptic vertex-weighted graphs and vertex-weighted graphs of genus at most $3$; see Remark~\ref{rmk:last:2}).  

Caporaso \cite[Conjecture~2.1]{Ca2} conjectured that 
the algebraic rank equals the combinatorial rank on 
any vertex-weighted graphs. It turns out that this is not the case in general; 
In \cite{CLM}, Caporaso, Len and Melo have found counterexamples,
which we have learned while preparing this article.
Since there are many graphs on which the algebraic rank equals the combinatorial rank (cf. Remark~\ref{rmk:last:2}), it will be an interesting question to characterize such graphs. 

\medskip
To prove Theorem~\ref{thm:cor:main}, we study the algebraic and combinatorial ranks of vertex-weighted graphs with a bridge. (For the definition of ${\rm Bs}(\left|\underline{d_i}\right|^\bullet)$, see Section~\ref{subsec:prelim:vw}.)  

\begin{Proposition}
\label{prop:main}
Let $\bar{G} = (G, \omega)$ be a vertex-weighted graph 
having a bridge $e$ with endpoints $v_1, v_2$.  
Let $G_1$ and $G_2$ be the connected components of
$G\setminus\{e\}$ 
such that 
$v_1 \in V(G_1), v_2\in V(G_2)$,  
and set $\bar{G}_i = (G_i, \rest{\omega}{V(G_i)})$ for $i = 1, 2$. 
Let $\underline{d} \in \Div(G)$, and let $\underline{d_i} \in \Div(G_i)$ 
be the restriction of $\underline{d}$ to $G_i$. Then we have 
\begin{equation}
\label{eqn:prop:main}
r_{\bar{G}}(\underline{d}) 
\leq 
\begin{cases}
r_{\bar{G}_1}(\underline{d_1}) +   r_{\bar{G}_2}(\underline{d_2}) +  1 
& \text{\textup{(}if $v_i \in {\rm Bs}(\left|\underline{d_i}\right|^\bullet)$ for each $i = 1, 2$\textup{)}, }\\ 
r_{\bar{G}_1}(\underline{d_1}) +   r_{\bar{G}_2}(\underline{d_2}) 
& \text{\textup{(}otherwise\textup{)}}. 
\end{cases}
\end{equation}
\end{Proposition}

There is a formula corresponding to \eqref{eqn:prop:main} (with the inequality replaced by the equality) for nodal curves (see Lemma~\ref{lemma:Inaba}). We prove Theorem~\ref{thm:cor:main} by 
induction on the number of bridges, 
using Proposition~\ref{prop:main}, Lemma~\ref{lemma:Inaba}  
and \cite[Corollary~1.7]{KY}. 

\medskip
To prove Theorem~\ref{thm:main:2}, we show 
the following proposition. (See Section~\ref{subsec:specialization:AC} for the notation.)  

\begin{Proposition}
\label{prop:KY:3}
Let $R$ be a complete discrete valuation ring 
with fractional field $\KK$ and residue field $k$. 
Let $\bar{G} = (G , \omega)$ be a non-hyperelliptic graph of 
genus $3$. Let $\Xscr$ be a regular, generically smooth, semi-stable $R$-curve with 
reduction graph $\bar{G}$.  Then the following condition 
\textup{(F)}  
holds. 
\begin{enumerate}
\item[(F)]
For any $\underline{d} \in \Div (G)$, there exists a divisor $\widetilde{D} \in \Div (\mathscr{X}_{\KK})$
such that $\widetilde{\rho}_*(\widetilde{D}) = \underline{d}$
and $r_{\bar{G}}(\underline{d}) = r_{\mathscr{X}_{\KK}}(\widetilde{D})$, 
\end{enumerate}
where $\mathscr{X}_{\KK}$ is the generic fiber of $\mathscr{X}$  
and $\widetilde{\rho}_*: \Div(\mathscr{X}_{\KK}) \to \Div(G)$ is the specialization map defined in \eqref{eqn:sp:2}.  
\end{Proposition}

We remark that a similar result for hyperelliptic graphs under a necessary  assumption on their bridges is obtained in \cite[Theorem~8.2]{KY}. (See Remark~\ref{rmk:last}. See also  
Proposition~\ref{prop:KY:3:2} for a related result, 
which says that any non-hyperelliptic graph of 
genus $3$ satisfies the condition (C) in \cite{KY}.) 
The proof of Proposition~\ref{prop:KY:3} uses the specialization lemma of Amini--Caporaso \cite[Theorem~4.10]{AC}, which is based on Baker's specialization lemma \cite{B}, and Raynaud's theorem on the surjectivity of 
the specialization map between principal divisors (see \cite{Ra}, 
\cite[Corollary~A2]{B} and Theorem~\ref{thm:Raynaud}).  
Then we deduce Theorem~\ref{thm:main:2} from 
Proposition~\ref{prop:KY:3} by the same argument 
as that in \cite{KY}, which is due to Caporaso.

\medskip
{\sl Acknowledgment.}\quad 
We deeply thank Professor Lucia Caporaso for 
invaluable comments on this article and \cite{KY}.  

%%%%%%%%%%%%%%%%%%%%%%%%%%%%%%%%%%%%%%%
%  Combinatorial and algebraic ranks of divisors on graphs    %
%%%%%%%%%%%%%%%%%%%%%%%%%%%%%%%%%%%%%
%
\setcounter{equation}{0}
\section{Combinatorial and algebraic ranks of divisors on graphs}
\label{sec:prelim}
In this section, we recall definitions and properties of combinatorial and algebraic ranks of divisors on graphs, 
which will be used later. 

\subsection{Divisors on finite graphs}
We briefly recall the theory of divisors on finite graphs. Our basic references are  
\cite{BN} and \cite{BN2}. 

Throughout this paper, a {\em finite graph} means an unweighted, finite connected graph. We allow a finite graph to have loops and multiple edges. For a finite graph $G$, let $V(G)$ denote the set of vertices, and $E(G)$ the set of edges. The {\em genus} of $G$ 
is defined as $g(G) = |E(G)| - |V(G)| + 1$. An edge $e \in E(G)$ 
is called a {\em bridge} if the deletion of $e$ makes $G$ disconnected. 

Let $\Div(G)$ be the free abelian group generated by $V(G)$. 
We call the elements of $\Div(G)$ {\em divisors} on $G$. 
Any divisor $\underline{d} \in \Div(G)$ is uniquely written as 
$\underline{d} = \sum_{v \in V(G)} n_v [v]$ for $n_v \in \ZZ$. 
The coefficient $n_v$ at $[v]$ is denoted by $\underline{d}(v)$. 
A divisor $\underline{d}$ is {\em effective}, 
written as $\underline{d} \geq 0$, 
if $\underline{d}(v) \geq 0$ for any $v \in V(G)$. 
The {\em degree} of a divisor $\underline{d}$ is defined as  
$\deg(\underline{d}) = \sum_{v \in V(G)} \underline{d}(v)$. 

A {\em rational function} on $G$ is an integer-valued function 
on $V(G)$. We denote by $\Rat(G)$ the set of rational functions on $G$. For 
$f \in \Rat(G)$ and a vertex $v$ of $G$, we set 
$
  \ord_v(f) = \sum_{e = \overline{wv} \in E(G)} \left(
  f(w) - f(v)
  \right), 
$
where the $e$'s run through all the edges of $G$ with endpoint $v$. Then 
\[
  \zero(f) := \sum_{v \in V(G)} \ord_v(f) [v] 
\]
is a divisor on $G$. 
The set of {\em principal divisors} on $G$ is defined as   
$ \Prin(G) = \{\zero(f) \mid f \in \Rat(G)\}$. Then $\Prin(G)$ is 
a subgroup of $\Div(G)$, and we write $\Pic(G) = \Div(G)/\Prin(G)$. 
For a divisor $\underline{d} \in \Div(G)$, 
let ${\rm cl}(\underline{d})$ denote
its divisor class in $\Pic(G)$.
% is usually 
%denoted by $\delta$, and we write $\underline{d} \in \delta$ and 
%$\delta = {\rm cl}(\underline{d})$. 

Two divisors $\underline{d} , \underline{d} ^\prime \in \Div(G)$ are said to be {\em linearly equivalent}, 
expressed as $\underline{d}  \sim \underline{d}^\prime$, if $\underline{d} - \underline{d}^\prime \in \Prin(G)$.  
For $\underline{d}  \in \Div(G)$, the complete linear system $|\underline{d}|$ is defined by 
\[
  |\underline{d} | = \{\underline{d}^\prime \in \Div(G) \mid \underline{d}^\prime  \geq 0, \quad \underline{d}^\prime  \sim \underline{d} \}. 
\]

\begin{Definition}[(Combinatorial) rank of a divisor \cite{BN}]
\label{def:rank:BN}
Let $G$ be a finite graph. 
Let $\underline{d} \in \Div(G)$. If $|\underline{d}| = \emptyset$, then we set 
$r_G(\underline{d}) := -1$. If $|\underline{d}| \neq \emptyset$, we set 
\[
 r_G(\underline{d}) := \max\left\{s \in \ZZ_{\geq 0} \;\left\vert\; 
\text{$|\underline{d}-\underline{e}| \neq \emptyset$ for 
any effective divisor $\underline{e}$ with $\deg(\underline{e}) = s$}
 \right. \right\}. 
\]
\end{Definition}

We note that $r_G(\underline{d})$ depends only on the divisor class of 
$\underline{d}$. For $\delta = {\rm cl}(\underline{d}) \in \Pic(G)$, we set 
$
  r_G(\delta) :=  r_G(\underline{d}).
$ 

A vertex $v \in V(G)$ is called a {\em base-point} of the complete linear system $|\underline{d} |$ if 
$r_G(\underline{d} - [v]) = r_G(\underline{d})$. The set of base-points of 
$|\underline{d} |$ is denoted by $\Bs(|\underline{d} |)$. 
If $|\underline{d} | = \emptyset$, then any vertex of $G$ is a base-point of $|\underline{d} |$
by definition. 

\medskip
{\em In the rest of this subsection, we assume that $G$ is loopless.} 
For any subset $A \subseteq V(G)$ and $v \in V(G)$, the {\em out-degree} of $v$ from $A$, denoted by $\mathrm{outdeg}_A(v)$, is the number of edges of $G$ having $v$ as one endpoint and whose other endpoint lies in $V(G)\setminus A$. 
For $\underline{d} \in \Div(G)$, a vertex $v \in A$ is  {\em saturated} 
for $\underline{d}$ with respect to $A$ if $\underline{d}(v) \geq \mathrm{outdeg}_A(v)$, and {\em non-saturated} otherwise. 

\begin{Definition}[$v_0$-reduced divisor \cite{BN}]
\label{def:reduced}
Fix a base vertex $v_0\in V(G)$.  
A divisor $\underline{d} \in \Div(G)$ is called a {\em $v_0$-reduced divisor} if 
$\underline{d}(v) \geq 0$ for any $v \in V(G)\setminus \{v_0\}$, and every non-empty subset $A$ of $V(G)\setminus \{v_0\}$ contains a non-saturated vertex 
$v \in A$ for $\underline{d}$ with respect to $A$. 
\end{Definition}

We recall from  \cite{BN} key properties of $v_0$-reduced divisors, which will be used later. 
 
\begin{Proposition}[{\cite[Proposition~3.1~and~its~proof]{BN}}] 
\label{prop:BN:reduced}
Fix a base vertex $v_0\in V(G)$.  Then for any 
$\underline{d} \in \Div(G)$, there exists a unique 
$v_0$-reduced divisor $\underline{d}^\prime \in \Div(G)$ that is 
linearly equivalent to $\underline{d}$. Further, $r_G(\underline{d}) \geq 0$ 
if and only if $\underline{d}^\prime$ is effective. 
\end{Proposition}

\medskip
The canonical divisor $K_G$ on $G$ is defined by 
$K_G = \sum_{v \in V(G)} (\val(v) - 2)[v] \in \Div(G)$, where $\val(v)$ 
denotes the number of edges with endpoint $v$. 
We remark that, with the above definition of rank, the notion of $v_0$-reduced divisors and the canonical divisor on $G$, Baker and Norine \cite[Theorem~1.12]{BN} established the Riemann--Roch theorem on a loopless finite graph.

Finally, we recall the definition of hyperelliptic graphs. 

\begin{Definition}[Hyperelliptic graph \cite{BN2}]
\label{def:hyp}
A loopless finite graph $G$ of $g(G) \geq 2$ is said to be {\em hyperelliptic} if 
there exists a divisor $\underline{d} \in \Div(G)$ such that 
$\deg(\underline{d}) = 2$ and $r_{\bar{G}}(\underline{d}) = 1$.  
\end{Definition}

%%%%%%%%%%%%%%%
\subsection{Rank of divisors on vertex-weighted graphs}
\label{subsec:prelim:vw}
We briefly recall the theory of divisors on vertex-weighted graphs. 
Our basic references are  
\cite{AC} and \cite{Ca2}. 

A {\em vertex-weighted} graph $\bar{G} = (G, \omega)$ is the pair of a finite graph $G$ and a function (called a vertex-weight function) $\omega: V(G) \to \ZZ_{\geq 0}$. The genus of $\bar{G}$ is defined as   
$g(\bar{G}) =  g(G) + \sum_{v \in V(G)} \omega(v)$. 

For a vertex-weighted graph $\bar{G} = (G, \omega)$, 
we make a loopless finite graph 
$\bar{G}^{\bullet}$ as follows: We add $\omega(v)$ loops 
to $G$ at $v$ for every vertex $v \in V(G)$; Then we 
insert a vertex in every loop edge. 
The graph $\bar{G}^{\bullet}$ is called {\em the 
virtual loopless finite graph of $\bar{G}$}. 

We have natural embeddings of the vertices 
$V(G) \subseteq V(\bar{G}^{\bullet})$, and of the 
divisor groups $\Div(G) \subseteq   \Div(\bar{G}^{\bullet})$. 
For $\underline{d} \in \Div(G)$, the rank $r_{\bar{G}}(\underline{d})$ of $\underline{d}$ is defined by 
\[
  r_{\bar{G}}(\underline{d}) := 
  r_{\bar{G}^{\bullet}}\left(\underline{d}\right), 
\]
where the right-hand side is defined in Definition~\ref{def:rank:BN}. 
Since $\Prin(G) \subseteq \Prin(\bar{G}^{\bullet})$, 
$r_{\bar{G}}(\underline{d})$ depends only on the divisor class of $\underline{d}$. 
For $\delta = {\rm cl}(\underline{d}) \in \Pic(G)$, 
we set 
$
  r_{\bar{G}}(\delta) := r_{\bar{G}}(\underline{d}). 
$

For $\underline{d} \in \Div(G)$, we write 
$|\underline{d}|^\bullet$ for the complete linear system $\left|\underline{d}\right|$ on $\bar{G}^{\bullet}$. Namely, we have 
\[
|\underline{d}|^\bullet := \{\underline{d}^\prime \in \Div(\bar{G}^{\bullet}) \mid 
\underline{d}^\prime \geq 0, \; \text{$\underline{d}^\prime$ is linearly equivalent to 
$\underline{d}$ in $\bar{G}^{\bullet}$}\}. 
\]
Here we use the notation ``$\bullet$'' to emphasize that we are considering divisors on $\bar{G}^{\bullet}$. 

Let $K_{\bar{G}^{\bullet}}$ be the canonical divisor of $\bar{G}^{\bullet}$. Then the support of $K_{\bar{G}^{\bullet}}$ lies in $V(G)$. We regard $K_{\bar{G}^{\bullet}}$ as an element of $\Div(G)$, and we define the {\em canonical divisor} 
$K_{\bar{G}}$ of $\bar{G}$ by $K_{\bar{G}} := K_{\bar{G}^{\bullet}} \in \Div(G)$. 
We remark that, if $G$ is loopless, then 
$K_{\bar{G}} 
= K_G + \sum_{v \in V(G)} 2 \omega(v) [v]$.

A vertex-weighted graph $\bar{G}$ of $g(\bar{G}) \geq 2$ is said to be {\em hyperelliptic} if its virtual loopless finite graph $\bar{G}^{\bullet}$ is hyperelliptic.

%%%%%%%%%%%%%%%
\subsection{Algebraic rank}
Following \cite{Ca2}, we recall the notion of 
the algebraic rank of a divisor class 
$\delta$ on a vertex-weighted graph. 

%As in the introduction, l
Let $k$ be a fixed algebraically closed field. 
By a {\em nodal curve}, we mean a connected, reduced, 
projective, one dimensional scheme over $k$ 
with at most ordinary double points as singularities. 

For a nodal curve $X$, the group of Cartier divisors is 
denoted by $\Div(X)$. We set $\Pic(X) = \Div(X)/\Prin(X)$, where 
$\Prin(X)$ denotes the group of principal divisors. For $L \in \Pic(X)$, 
we write $r_X(L) = \dim_k H^0(X, L) - 1$. 

Given a nodal curve $X$, 
the (vertex-weighted) {\em dual graph} $\bar{G} = (G, \omega)$ associated to 
$X$ is defined as follows.  
Let $C_1, \ldots, C_r$ be the irreducible components of $X$. Then 
$G$ has vertices $v_1, \ldots, v_r$ which correspond to $C_1, \ldots, C_r$, respectively. Two vertices $v_i, v_j$ ($i \neq j$) of $G$ are connected by 
$a_{ij}$ edges if $\# C_i \cap C_j = a_{ij}$.  A vertex $v_i$ has $b_i$ loops if 
$\#{\rm Sing}(C_i) = b_i$.  The vertex-weighted function $\omega$ is given 
by assigning to $v_i$ the geometric genus of $X_i$. 

Let $\bar{G} = (G, \omega)$ be a vertex-weighted graph. 
Let $M^{\alg}(\bar{G})$ be a family of nodal curves 
representing all the isomorphism classes of
nodal curves with
dual graph $\bar{G}$.
For $X \in M^{\alg}(\bar{G})$, we write  $X = \cup_{v \in V(G)} C_v$, where $C_v$ is the irreducible curve corresponding to $v \in V(G)$. We have a natural map 
\begin{equation}
\label{eqn:sp}
  \rho_*: \Div(X) \to \Div(G), \quad
  D \mapsto \sum_{v \in V(G)} \left(\deg\left(\rest{D}{C_v}\right)\right)[v]. 
\end{equation}
In other words, for a Cartier divisor $D$ on $X$, 
$\rho_*(D) \in \Div(G)$ gives the multidegree of $D$. 
Since linear equivalent divisors on $X$ have the same 
multidegree, $\rho_*$ descends to $\Pic(X) \to \Div(G)$. Then we have a stratification of $\Pic(X)$: 
\[
  \Pic(X) = \bigsqcup_{\underline{d} \in \Div(G)} \Pic^{\underline{d}}(X),  
\]
where $\Pic^{\underline{d}}(X) 
= \left\{L \in \Pic(X) \mid \text{$\deg\left(\rest{L}{C_v}\right) = d_v$ for any 
$v \in V(G)$} \right\}$ for $\underline{d} = (d_v)_{v \in V(G)} \in  \Div(G)$. 

\begin{Definition}[Algebraic rank \cite{Ca2}]
Let $\bar{G} = (G, \omega)$ be a vertex-weighted graph, and 
$\delta \in \Pic(G)$ a divisor class on $G$. We set 
\[
r_{\bar{G}}^{\alg}(\delta)
=  \max_{X \in M^{\alg}(\bar{G})} \left\{
  \min_{\underline{d} \in \delta} \left\{
  \max_{L \in \Pic^{\underline{d}}(X)} \left\{
  r_X(L)
  \right\}
  \right\}
  \right\}, 
\]
and call $r_{\bar{G}}^{\alg}(\delta)$ the {\em algebraic rank} 
of combinatorial type $\delta$. 
\end{Definition}

%%%%%%%%%%%%%%%%%%%%%%
\subsection{The specialization lemma for vertex-weighted graphs}
\label{subsec:specialization:AC}
We recall the specialization lemma for vertex-weighted graphs 
due to Amini--Caporaso \cite{AC}, which generalizes 
Baker's specialization lemma for loopless finite graphs \cite{B}. 
Our basic references are \cite{AC} and \cite{B}. 

Let $k$ be a fixed algebraically closed field. 
Let $R$ be a complete discrete valuation ring with residue field $k$. 
Let $\KK$ denote the fractional field of $R$. 

By an {\em $R$-curve}, we mean an integral scheme of dimension 
$2$ that is proper and flat over $\Spec(R)$. For an $R$-curve $\Xscr$, 
we denote by $\Xscr_{\KK}$ the generic fiber of $\Xscr$, and 
by $X$ the special fiber of $\Xscr$. We say that $\Xscr$ is a {\em semi-stable} $R$-curve if $X$ is a nodal curve. The vertex-weighted dual graph $\bar{G} = (G, \omega)$ of $X$ is then called the {\em reduction graph} of $\Xscr$. 

Let $\Xscr$ be a regular, generically smooth, semi-stable $R$-curve. Since $\Xscr_{\KK}$ is smooth (resp. $\Xscr$ is regular), the group of Cartier divisors 
on $\Xscr_{\KK}$ (resp. $\Xscr$) is the same as the group of Weil divisors. 
The Zariski closure of an effective divisor on $\Xscr_{\KK}$ in $\Xscr$ is a 
Cartier divisor. Extending by linearity, one can associate 
to any divisor on $\Xscr_{\KK}$ a Cartier divisor on $\Xscr$, which is also called the {Zariski closure} of the divisor. 

Let $\widetilde{D}$ be a divisor on $\Xscr_{\KK}$ 
and $\widetilde{\Dscr}$ the Zariski closure 
of $\widetilde{D}$. Let $\OO_{\Xscr}(\widetilde{\Dscr})$ be the invertible sheaf on $\Xscr$ associated to $\widetilde{\Dscr}$. 
We define the {\em specialization map} 
$\widetilde{\rho}_*: \Div(\Xscr_{\KK}) \to \Div(G)$ by 
\begin{equation}
\label{eqn:sp:2}
  \widetilde{\rho}_*(\widetilde{D}) :=  
  \sum_{v \in V(G)} \deg\left(
  \OO_{\Xscr}(\widetilde{\Dscr})\vert_{C_v}
  \right) [v]
  \in \Div(G) 
\end{equation}
(see \cite[\S2.1]{B}). The map $\widetilde{\rho}_*$ is compatible with 
the map $\rho_*$ in \eqref{eqn:sp}: Namely, let $D \in \Div(X)$ be a Cartier 
divisor on the special fiber such that the associated invertible sheaf $\OO_X(D)$ 
is isomorphic to $\rest{\OO_{\Xscr}(\widetilde{\Dscr})}{X}$; Then, by definition, we have 
\begin{equation}
\label{eqn:compatibility}
  \rho_*(D) = \widetilde{\rho}_*(\widetilde{D}).
\end{equation}  

\begin{Remark}
\label{rmk:notation}
In \cite{KY}, $\widetilde{\rho}_*$ is denoted by $\rho_*$. Here we use the notation 
$\widetilde{\rho}_*$, for we have already use the notation $\rho_*$ in \eqref{eqn:sp}. 
\end{Remark}

Let $\Div(\Xscr(\KK))$ be the subgroup 
of $\Div(\Xscr_{\KK})$ generated by 
$\KK$-valued points of $\Xscr$. 
Then 
\begin{equation}
\label{eqn:surj:rho}
  \rest{\widetilde{\rho}_*}{\Div(\Xscr_{\KK})}: 
  \Div(\Xscr(\KK)) \to \Div(G)
\end{equation}
is surjective (see \cite[Remark~2.3]{B} and \cite[Proposition~10.1.40(b)]{Liu}). 

\begin{Theorem}[Amini--Caporaso's~specialization~lemma~{\cite[Theorem~4.10]{AC}}]
\label{thm:AC}
Let $\Xscr$ be a regular, generically smooth, semi-stable $R$-curve 
with reduction graph $\bar{G} = (G, \omega)$. 
Then, for any $\widetilde{D} \in \Div(\Xscr_\KK)$, one has 
$r_{\bar{G}}(\widetilde{\rho}_*(\widetilde{D})) 
\geq r_{\Xscr_\KK}(\widetilde{D})$. 
\end{Theorem}

Theorem~\ref{thm:AC} is a generalization of Baker's specialization lemma for loopless finite graphs in \cite{B}.  
Although Amini and Caporaso consider 
a smooth quasi-projective curve $B$ over $k$ (in place of $\Spec(R)$), 
i.e., they consider a morphism $\phi: \Xscr \to B$, 
we remark that their arguments 
also work over $\Spec(R)$. (By the surjectivity of the map \eqref{eqn:surj:rho}, the argument over $R$ works 
as the same as the argument for $\phi: \Xscr \to B$ 
which admits a section passing through any given component 
of the special fiber.)

%%%%%%%%%%%%%%%%%%%%%%%
%  Proof of Proposition~\ref{prop:main}    %
%%%%%%%%%%%%%%%%%%%%%%%
\setcounter{equation}{0}
\section{Reduced divisors and decomposition of graphs}
In this section, we prove Proposition~\ref{prop:main}. 
We first show some properties of divisors on a graph with a bridge. 

\begin{Lemma}
\label{lemma:red:divisors}
Let $G$ be a loopless finite graph 
with a bridge $e$ having endpoints $v_1, v_2$. 
Let $G_1$ and $G_2$ be the connected components of
$G\setminus\{e\}$ such that $v_1 \in V(G_1), v_2\in V(G_2)$. 
For $i = 1, 2$, 
let $\jmath_i: V(G_i) \hookrightarrow V(G)$ be the natural embedding 
and $\jmath_{i*}: \Div(G_i) \hookrightarrow \Div(G)$ the induced map. 
\begin{enumerate}
\item
For $i = 1, 2$, 
we have $\jmath_{i*}\left(\Prin(G_i)\right) \subseteq \Prin(G)$. 
\item
For $i = 1, 2$, 
let $\underline{d_i}$ be a $v_i$-reduced divisor 
on $G_i$.  Then 
\begin{equation}
\label{eqn:red:divisor}
  \jmath_{1*}\left(\underline{d_1} - \underline{d_1}(v_1)[v_1]\right) 
  + \jmath_{2*}\left(\underline{d_2} - \underline{d_2}(v_2)[v_2]\right) 
  + \left(\underline{d_1}(v_1) + \underline{d_2}(v_2)\right)[v_1]
\end{equation}
is a $v_1$-reduced divisor on $G$.  
\end{enumerate}
\end{Lemma}

\Proof
(1) We may assume that $i = 1$. Let $f_1$ be a rational function 
on $G_1$. We extend $f_1$ to a rational function $\widetilde{f}_1$ on $G$ 
by setting $\widetilde{f}_1(w) = f_1(v_1)$ for any $w \in V(G_2)$. Then we have 
$
  \zero(\widetilde{f}) = j_{1*}\left(
  \zero(f_1)
  \right)
$. Thus $j_{1*}\left(
  \zero(f_1)
  \right) \in \Prin(G)$, which gives the assertion. 

(2) 
We put
\[
\underline{d} :=
\jmath_{1*}\left(\underline{d_1} - \underline{d_1}(v_1)[v_1]\right) 
  + \jmath_{2*}\left(\underline{d_2} - \underline{d_2}(v_2)[v_2]\right) .
\]
It suffices to show that $\underline{d}$ is
a
$v_1$-reduced divisor on $G$.  
%We denote the divisor in \eqref{eqn:red:divisor} by 
%$\underline{d}$. 
%To show this, 
Let $A \subseteq V(G)\setminus\{v_1\}$ be any non-empty 
subset,
and we are going 
show that there exists a non-saturated vertex 
$v \in A$ for $\underline{d}$ with respect to $A$. 

If $v_2 \in A$, then it follows from $v_1 \not\in A$ 
that $\outdeg_{A}(v_2) \geq 1$ (from the contribution of the bridge $e$). 
Since $\underline{d}(v_2) = 0$, we see that 
$v_2 \in V(G) \setminus \{v_1\}$ is 
a non-saturated vertex for $\underline{d}$ with respect to $A$. 
Thus we may and do assume  that $v_2 \not\in A$,
and hence $A \subseteq V(G)\setminus\{v_1 , v_2\}$.

We set $A_1 := A \cap V(G_1)$ and $A_2 := A \cap V(G_2)$. Then 
$A_1 \subseteq V(G)\setminus\{v_1\}$ and 
$A_2 \subseteq V(G)\setminus\{v_2\}$. 
Since $A \neq \emptyset$, we have  
$A_1 \neq \emptyset$ or $A_2 \neq \emptyset$.  
Without loss of generality,
we assume that
%Suppose that 
$A_1 \neq \emptyset$. 
Since $\underline{d_1}$ is a $v_1$-reduced divisor on 
$G_1$, there exists a non-saturated vertex $v \in A_1$ 
for $\underline{d_1}$ with respect to $A_1$, i.e., 
$\underline{d_1}(v) < \outdeg_{A_1}(v)$. Since 
$\underline{d_1}(v) = \underline{d}(v)$ and 
$\outdeg_{A_1}(v) = \outdeg_{A}(v)$, 
we have $\underline{d}(v) < \outdeg_{A}(v)$. 
Thus $v \in A$ is a non-saturated vertex 
for $\underline{d}$ with respect to $A$,
which shows the lemma.
%Suppose that $A_2 \neq \emptyset$. 
%Then, by a similar argument as above, 
%we see that there exists a non-saturated vertex $v \in A$ is a non-saturated vertex 
%for $\underline{d}$ with respect to $A$. 
%This 
%concludes that $\underline{d}$ is a $v_1$-reduced divisor on $G$. 
\QED

The next lemma will be used in Section~\ref{sec:genus:3}. 

\begin{Lemma}
\label{lemma:for:genus:3}
Let $\bar{G} = (G, \omega)$ be a vertex-weighted graph. 
Let $\underline{d} \in \Div(G)$. If $r_{\bar{G}}(\underline{d}) \geq 0$, 
then there exists an effective divisor $\underline{e} \in \Div(G)$ 
that is linearly equivalent to $\underline{d}$ in $G$. 
\end{Lemma}

\Proof
Let $\bar{G}^\bullet$ be the virtual loopless finite graph of $\bar{G}$. 
Via the natural embedding of the sets of vertices, we regard 
$V(G) \subseteq V(\bar{G}^\bullet)$. The condition 
$r_{\bar{G}}(\underline{d}) := r_{\bar{G}^\bullet}(\underline{d}) \geq 0$ means that there exists a rational function $\bar{f} \in \Rat(\bar{G}^\bullet)$ such that $\underline{d}^\prime := \underline{d} + \zero(\bar{f})$ is an effective divisor on $\bar{G}^\bullet$. 

Let $w \in V(\bar{G}^\bullet)\setminus V(G)$. This means that $w$ is a vertex inserted in a loop edge.  Thus there exist exactly two edges $e_1, e_2$ of $\bar{G}^\bullet$ with endpoint $w$, and the other endpoint of $e_1$ and that of $e_2$ are the same, which we denote by $w^\prime$. Since $\underline{d}(w) = 0$ and 
$\underline{d}^\prime(w) \geq 0$, we see that 
$\bar{f}(w^\prime) \geq \bar{f}(w)$.
% for every $w \in V(\bar{G}^\bullet)\setminus V(G)$. 

We set $f := \rest{\bar{f}}{V(G)} \in \Rat(G)$. Since $\underline{d} \in \Div(G)$ and $\bar{f}(w^\prime) \geq \bar{f}(w)$ for every $w \in V(\bar{G}^\bullet)\setminus V(G)$, we see that 
$
 \underline{e} := \underline{d} + \zero(f)  
$
is an effective divisor on $G$. 
This shows the lemma.
\QED

We begin the proof of Proposition~\ref{prop:main}. 

\medskip
{\bf Proof of Proposition~\ref{prop:main}.}\quad
Let $\bar{G}^\bullet, \bar{G}_1^\bullet$ and $\bar{G}_2^\bullet$ be the virtual loopless finite graphs of $\bar{G}, \bar{G}_1$ and $\bar{G}_2$, respectively. 
Note that $\bar{G}^\bullet$ is the graph 
obtained by connecting $\bar{G}_1^\bullet$ and $\bar{G}_2^\bullet$
with the edge $e$. 
For $i = 1, 2$, 
let $\jmath_{i*}^\bullet: \Div(\bar{G}_i^\bullet) \hookrightarrow 
\Div(\bar{G}^\bullet)$ be the induced embedding of divisors.

Via the natural embedding of the sets of vertices, we regard 
$V(G_i)\subseteq  V(G) \subseteq V(\bar{G}^\bullet)$ and 
$V(G_i)\subseteq V(\bar{G}_i^\bullet) \subseteq V(\bar{G}^\bullet)$
for $i = 1, 2$.
Thus,
in the following argument,
we
will often identify the vertex $v_i \in V(G_i)$ with 
the corresponding vertices in $G$, $\bar{G}_i^\bullet$ and $\bar{G}^\bullet$.

%Notice that $\bar{G}^\bullet$ is the graph that connects 
%$\bar{G}_1^\bullet$ and $\bar{G}_2^\bullet$ by the edge $e$. 
%For $i = 1, 2$, we write 
%$\jmath_i^\bullet: V(\bar{G}_i^\bullet) \hookrightarrow 
%V(\bar{G}^\bullet)$ for the natural embedding of the sets of vertices. 
%Let $\jmath_{i*}^\bullet: \Div(\bar{G}_i^\bullet) \hookrightarrow 
%\Div(\bar{G}^\bullet)$ be the induced embedding of divisors. 

For $i = 1, 2$, 
we set $r_i = r_{\bar{G}_i}(\underline{d_i})$. 
By the definition of the rank, there exists an effective 
divisor $\underline{e_i} \in \Div(G_i^\bullet)$ with 
$\deg(\underline{e_i}) = r_i + 1$ such that 
$r_{\bar{G}_i^\bullet}(\underline{d_i} - \underline{e_i}) = -1$.  
We set $\underline{f_i} =\underline{d_i} - \underline{e_i}$, and 
let $\underline{f_i}^\prime \in \Div(G_i^\bullet)$ be the $v_i$-reduced 
divisor that is linearly equivalent to $\underline{f_i}$ on $G_i^\bullet$. 
Since $r_{\bar{G}_i^\bullet}(\underline{d_i} - \underline{e_i}) = -1$, 
we have $\underline{f_i}^\prime(v_i) <0$ by Proposition~\ref{prop:BN:reduced}. 

%We are going to show 
We claim that 
$ r_{\bar{G}^\bullet}\left(\underline{d} - \jmath_{1*}^\bullet(\underline{e_1})- \jmath_{2*}^\bullet(\underline{e_2})\right) = -1$. 
Indeed,
we see from Lemma~\ref{lemma:red:divisors}(1) that, as divisors 
on $\bar{G}^\bullet$, 
{\allowdisplaybreaks
\begin{align*}
& \underline{d} 
- \jmath_{1*}^\bullet(\underline{e_1})
- \jmath_{2*}^\bullet(\underline{e_2}) \\
&\qquad =  
 \jmath_{1*}^\bullet\left(\underline{d_1} - \underline{e_1} \right) 
+  \jmath_{2*}^\bullet\left(\underline{d_2}- \underline{e_2}\right) \\
&  \qquad\sim 
 \jmath_{1*}^\bullet\left(\underline{f_1}^\prime\right) 
 + \jmath_{2*}^\bullet\left(\underline{f_2}^\prime\right) \\
&  \qquad= 
\jmath_{1*}^\bullet\left(\underline{f_1}^\prime 
- \underline{f_1}^\prime(v_1)[v_1] \right) 
+ \jmath_{2*}^\bullet\left(\underline{f_2}^\prime 
- \underline{f_2}^\prime(v_2)[v_2] \right)
+ (\underline{f_1}^\prime(v_1)+ \underline{f_2}^\prime(v_2))[v_1] .
%\\
%&  \qquad=:  \underline{g}.  
\end{align*}
}
We denote by $\underline{g}$ the divisor 
in the last line in the above. 
Then, by
Lemma~\ref{lemma:red:divisors}(2), 
$\underline{g}$ is a $v_1$-reduced divisor on $\bar{G}^\bullet$. 
Since $\underline{g}(v_1) = \underline{f_1}^\prime(v_1) + \underline{f_2}^\prime(v_2) <0$, we have $ r_{\bar{G}^\bullet}\left(\underline{d} - \jmath_{1*}^\bullet(\underline{e_1})- \jmath_{2*}^\bullet(\underline{e_2})\right) = -1$  by Proposition~\ref{prop:BN:reduced}. 

It follows that 
\begin{align*}
  r_{\bar{G}}(\underline{d}) 
  & = 
  r_{\bar{G}^\bullet}\left(\underline{d}\right) 
  \leq \deg\left(\jmath_{1*}^\bullet(\underline{e_1}) + \jmath_{2*}^\bullet(\underline{e_2})\right) - 1 \\
  & = 
  r_1 + r_2 + 1 
  = r_{\bar{G}_1}(\underline{d_1}) 
  + r_{\bar{G}_2}(\underline{d_2}) + 1. 
\end{align*}
This shows the inequality in 
(\ref{eqn:prop:main})
in the case of  
$v_i \in {\rm Bs}(\left|\underline{d_i}\right|^\bullet)$ for each $i = 1, 2$.

Suppose now that $v_1 \not\in {\rm Bs}(|\underline{d}_1|^\bullet)$ or 
$v_2 \not\in {\rm Bs}(|\underline{d}_2|^\bullet)$. We need to show that 
\[
  r_{\bar{G}}(\underline{d}) \leq r_{\bar{G}_1}(\underline{d_1}) 
  + r_{\bar{G}_2}(\underline{d_2}). 
\]
Without loss of generality, 
we may assume that $v_2 \not\in {\rm Bs}(|\underline{d}_2|^\bullet)$. 
This means that 
%$|\underline{d}_2|^\bullet \neq \emptyset$ on 
%$\bar{G}_2^\bullet$ and 
$r_{\bar{G}_2}(\underline{d_2} - [v_2]) 
= r_2 - 1$. 
Note that $r_2 \geq 0$.

By the definition of the rank, there exists an effective 
divisor $\underline{\widetilde{e}_2} \in \Div(\bar{G}_2^\bullet)$ with 
$\deg(\underline{\widetilde{e}_2}) = r_2$ such that 
$ r_{\bar{G}_2^\bullet}\left(\underline{d_2} - [v_2] - \underline{\widetilde{e}_2}\right) = -1$.  
We set $\underline{h_2} = \underline{d_2} - \underline{\widetilde{e}_2}$, and 
let $\underline{h_2}^\prime \in \Div(\bar{G}_2^\bullet)$ be the $v_2$-reduced 
divisor that is linearly equivalent to $\underline{h_2}$ on $\bar{G}_2^\bullet$. 
Since $r_{\bar{G}_2}(\underline{d_2}) = r_2$, 
we have $r_{\bar{G}_2^\bullet}(\underline{d_2} - \underline{\widetilde{e}_2}) \geq 0$, 
hence $\underline{h_2}^\prime$ is an effective divisor on $\bar{G}_2^\bullet$ by 
Proposition~\ref{prop:BN:reduced}. 
Since $\underline{h_2}^\prime$ is $v_2$-reduced 
and since $
r_{\bar{G}_2^\bullet}(\underline{h_2}^\prime - [v_2]) = 
r_{\bar{G}_2^\bullet}(\underline{d_2} - [v_2] - \underline{\widetilde{e}_2}) = -1$, 
we see that $\underline{h_2}^\prime(v_2) = 0$. 

%We are going to show that 
%$r_{\bar{G}^\bullet}\left(\underline{d} - \jmath_{1*}^\bullet(\underline{e_1})- \jmath_{2*}^\bullet(\widetilde{e}_2)\right) = -1$
It follows from Lemma~\ref{lemma:red:divisors}(1) that  
\begin{align*}
  \underline{d} - \jmath_{1*}^\bullet(\underline{e_1})- \jmath_{2*}^    \bullet(\widetilde{e}_2)
 &  = 
  \jmath_{1*}^\bullet\left(\underline{d_1} - \underline{e_1} \right) 
+  \jmath_{2*}^\bullet\left(\underline{d_2}- \underline{\widetilde{e}_2}\right)   \\
 &  \sim
  \jmath_{1*}^\bullet\left(\underline{f_1}^\prime\right)   + 
  \jmath_{2*}^\bullet\left(\underline{h_2}^\prime\right). 
\end{align*}
Since $\underline{h_2}^\prime(v_2) = 0$, we see that 
$\jmath_{1*}^\bullet\left(\underline{f_1}^\prime\right)   + 
  \jmath_{2*}^\bullet\left(\underline{h_2}^\prime\right)$ is a 
$v_1$-reduced divisor on $\bar{G}^\bullet$ 
by Lemma~\ref{lemma:red:divisors}(2). 
Since $
\underline{f_1}^\prime(v_1) + \underline{h_2}^\prime(v_2)
= \underline{f_1}^\prime(v_1) <0$, Proposition~\ref{prop:BN:reduced} tells us that 
$r_{\bar{G}^\bullet}\left(\underline{d} - \jmath_{1*}^\bullet(\underline{e_1})- \jmath_{2*}^\bullet(\widetilde{e}_2)\right) = -1$. 
Then 
\begin{align*}
  r_{\bar{G}}(\underline{d}) 
  & = r_{\bar{G}^\bullet}(\underline{d}) 
  \leq \deg( \jmath_{1*}^\bullet(\underline{e_1}) + \jmath_{2*}^\bullet(\widetilde{e}_2)) -1 \\
  & = r_1 + r_2 = r_{\bar{G}_1}(\underline{d_1}) 
  + r_{\bar{G}_2}(\underline{d_2}). 
\end{align*}
Thus we obtain the inequality
in the remaining case.
\QED

There exists 
a formula corresponding to Proposition~\ref{prop:main} (with the inequality replaced by the equality) for nodal curves. 

\begin{Lemma}
\label{lemma:Inaba}
Let $X$ be a nodal curve. 
We assume that $X$ has a decomposition as $X = X_1 \cup X_2$ into 
two nodal curves so that $X_1$ and $X_2$ meet
at exactly one point $p$. Let $D$ be a Cartier divisor on $X$, and we set 
$D_i = \rest{D}{X_i} \in \Div(X_i)$ for $i = 1, 2$.  Then 
\begin{equation}
\label{eqn:lemma:Inaba}
r_{X}(D) 
= 
\begin{cases}
r_{X_1}(D_1) +   r_{X_2}(D_2) +  1 
& \text{\textup{(}if $p \in {\rm Bs}(\left|D_i\right|)$ for each $i = 1, 2$\textup{)}, }\\ 
r_{X_1}(D_1) +   r_{X_2}(D_2) 
& \text{\textup{(}otherwise\textup{)}}. 
\end{cases}
\end{equation}
\end{Lemma}

\Proof
This is a well-known fact, so 
we omit a proof. See also \cite[Remark~1.5]{Ca2}. 
%For $i = 1, 2,$, let $\jmath_i: X_i \to X$ be the natural emedding. 
%Let $\OO_X(D)$ be the locally free sheaf on $X$ associated to 
%$D$, and $\OO_{X_i}(D_i)$ be the locally free sheaf on $X_i$ associated to 
%$D_i$. Then we have the natural identification 
%$\rest{\OO_X(D)}{p} \cong \rest{\OO_{X_1}(D_1)}{p} \cong 
%\rest{\OO_{X_2}(D_2)}{p}$, which we regard as 
%the skyscraper sheaf $k_p$ on $X$. 
%
%We have an exact sequence
%\[
% 0 \to \OO_X(D) \overset{\varphi}{\longrightarrow} \jmath_{1*}(\OO_{X_1}(D_1)) \oplus 
% \jmath_{2*}(\OO_{X_2}(D_2)) \overset{\psi}{\longrightarrow} k_p \to 0, 
%\]
%where, locally, $\varphi$ is given by $s \mapsto (\rest{s}{X_1}, \rest{s}{X_2})$ and 
%$\psi$ is given by $(t_1, t_2) \mapsto \rest{t_1}{p} -  \rest{t_2}{p}$. 
%We obtain the long exact sequence
%\[
% 0 \to H^0(X, \OO_X(D)) \to 
% H^0(X_1, \OO_{X_1}(D_1)) \oplus H^0(X_2, \OO_{X_2}(D_2))
% \to k \to \cdots 
%\]
%We note that the map $H^0(X_i, \OO_{X_i}(D_i)) \to k$ given by 
%$t \mapsto t(p)$ is zero map if and only if $p$ is a base-point of $|D_i|$. 
%Thus we get 
%\begin{multline*}
%h^0(X, \OO_X(D)) \\
%= 
%\begin{cases}
%h^0(X_1, \OO_{X_1}(D_1)) +   h^0(X_2, \OO_{X_1}(D_2)) 
%& \text{\textup{(}if $p \in {\rm Bs}(\left|D_i\right|)$ for each $i = 1, 2$\textup{)}, }\\ 
%h^0(X_1, \OO_{X_1}(D_1)) +   h^0(X_2, \OO_{X_1}(D_2)) - 1
%& \text{\textup{(}otherwise\textup{)}}. 
%\end{cases}
%\end{multline*}
%Since $r_X(D) = h^0(X, \OO_X(D)) - 1$ and 
%$r_{X_i}(D_i) = h^0(X_i, \OO_{X_i}(D_i)) - 1$, 
%we obtain \eqref{eqn:lemma:Inaba}. 
\QED

The following simple remark will be used in the next section. 

\begin{Remark}
\label{rmk:simple}
Let $X, X_1, X_2, p$ be as in Lemma~\ref{lemma:Inaba}. For $i = 1, 2$, 
let $D_i$ be a Cartier divisor on $X_i$. Then there exists a Cartier 
divisor $D$ on $X$ such that $\rest{D}{X_i}$ is linearly equivalent to $D_i$. 
Indeed, let $p_i: X \to X_i$ be the morphism given by the identity on $X_i$ 
and the constant map to $p$ on the other component. Let $\OO_{X_i}(D_i)$ be the 
invertible sheaf on $X_i$ associated to $D_i$. Then it suffices to take 
$D \in \Div(X)$ such that the associated invertible sheaf $\OO_X(D)$ is isomorphic to  $p_1^*\left(
\OO_{X_1}(D_1)
\right) \otimes p_2^*\left(
\OO_{X_2}(D_2)
\right)$. 
\end{Remark}

%%%%%%%%%%%%%%%%%%%%%%%%%%%
%  Graphs with a bridge and hyperelliptic graphs %
%%%%%%%%%%%%%%%%%%%%%%%%%%%
\setcounter{equation}{0}
\section{Graphs with a bridge and hyperelliptic graphs}
\label{sec:bridge}
In this section, we prove Theorem~\ref{thm:cor:main}. 
We begin by showing the 
following lemma. 

\begin{Lemma}
\label{lemma:useful}
Let $\bar{G} = (G, \omega)$ be a vertex-weighted graph 
with a bridge $e$ with endpoints $v_1, v_2$. Let 
$G_1$ and $G_2$ be the connected components of
$G\setminus\{e\}$ such that 
$v_1 \in V(G_1)$ and $v_2 \in V(G_2)$, and we set 
$\bar{G}_i = (G_i, \rest{\omega}{V(G_i)})$ for $i = 1, 2$. 
Let $\underline{d} \in \Div(G)$, and let $\underline{d_i} \in \Div(G_i)$ be 
the restriction of $\underline{d}$ to $G_i$. 
Let $X$ be a nodal curve over $k$ with dual graph $\bar{G}$, and 
we write $X_i$ for 
the union of irreducible components of $X$ corresponding to $\bar{G}_i$. 
Let $\rho_*: \Div(X) \to \Div(G)$ and 
$\rho_{i*}: \Div(X_i) \to \Div(G_i)$ be the maps defined in \eqref{eqn:sp}. 
For $i = 1, 2$, we 
assume that, for any divisor $\underline{e_i} \in \Div(G_i)$, 
there exists a Cartier divisor $E_i$ on $X_i$ 
satisfying $\rho_{i*}(E_i) = \underline{e_i}$ and 
$r_{X_i}(E_i) \geq r_{\bar{G}_i}(\underline{e_i})$. 
Then there exists a Cartier divisor $D$ on $X$ satisfying 
$\rho_{*}(D) = \underline{d}$ and 
\[
 r_X(D) \geq 
 \begin{cases}
 r_{\bar{G}_1}(\underline{d_1} - [v_1]) + 
 r_{\bar{G}_2}(\underline{d_2} - [v_2])  + 2 &  
 \text{\textup{(}if $v_i \not\in {\rm Bs}(\left|\underline{d_i}\right|^\bullet)$ for each $i = 1, 2$\textup{)}, } \\
 r_{\bar{G}_1}(\underline{d_1} - [v_1]) + 
 r_{\bar{G}_2}(\underline{d_2} - [v_2]) + 1 & 
 \text{\textup{(}otherwise\textup{)}}. 
 \end{cases}
\]
\end{Lemma}
 
\Proof
{\bf Case 1.}\quad 
Suppose that $v_i \not\in {\rm Bs}(\left|\underline{d_i}\right|^\bullet)$ for each $i = 1, 2$. 
This means that $\left|\underline{d_i}\right|^\bullet \neq \emptyset$ and 
\begin{equation}
\label{eqn:lem:usuful:1}
r_{\bar{G}_i}(\underline{d_i} - [v_i])
 = r_{\bar{G}_i}(\underline{d_i}) - 1. 
\end{equation}
We take a Cartier divisor $D_i$ on $X_i$ 
satisfying $\rho_{i*}(D_i) = \underline{d_i}$ and 
$r_{X_i}(D_i) \geq r_{\bar{G}_i}(\underline{d_i})$. 
By Remark~\ref{rmk:simple}, there exists a Cartier divisor $D$ on $X$ such that 
$\rest{D}{X_i}$ is linearly equivalent to $D_i$ on $X_i$. Then we have 
\begin{align*}
  \rho_{*}(D) 
  & = \rho_{1*}(\rest{D}{X_1}) + \rho_{2*}(\rest{D}{X_2}) \\
  & = \rho_{1*}(D_1) + \rho_{2*}(D_2)  
  =  \underline{d_1} + \underline{d_2}  = \underline{d}. 
\end{align*}
Further, we have 
\begin{align*}
  r_X(D) 
  & \geq r_{X_1}(\rest{D}{X_1}) + r_{X_2}(\rest{D}{X_2}) 
   && \textup{(from Lemma~\ref{lemma:Inaba})} \\
  & = r_{X_1}(D_1) + r_{X_2}(D_2) 
  &&  \textup{(since $\rest{D}{X_i} \sim D_i$ for each $i = 1, 2$)} \\
  & \geq r_{\bar{G}_1}(\underline{d_1}) + r_{\bar{G}_2}(\underline{d_2})
  && \textup{(from the assumptions on $D_i$)} \\
  & = r_{\bar{G}_1}(\underline{d_1} - [v_1]) + 
 r_{\bar{G}_2}(\underline{d_2} - [v_2])  + 2 && 
 \textup{(from \eqref{eqn:lem:usuful:1})}. 
\end{align*}
This gives the desired properties in this case. 

\smallskip
{\bf Case 2.}\quad 
Suppose that  
$v_1\in {\rm Bs}(\left|\underline{d_1}\right|^\bullet)$ or 
$v_2\in {\rm Bs}(\left|\underline{d_2}\right|^\bullet)$. 
For $ i = 1, 2$, 
we take a Cartier divisor $D_i^\prime$ on $X_i$ 
satisfying $\rho_{i*}(D_i^\prime) = \underline{d_i} - [v_i]$ and 
$r_{X_i}(D_i^\prime) \geq r_{\bar{G}_i}(\underline{d_i} - [v_i])$. 

We remark that $X = X_1 \cup X_2$
and that $X_1 \cap X_2$ consists of the node
of $X$ corresponding
to the edge $e$.
Let $p$ denote this node.
Since $p$ is a smooth point on $X_i$, 
the Weil divisor $[p]$ is regarded as a Cartier divisor on $X_i$. 
We set
\[
  D_i = D_i^\prime + [p] \in \Div(X_i). 
\]
By Remark~\ref{rmk:simple}, there exists a Cartier divisor $D$ on $X$ such that 
$\rest{D}{X_i}$ is linearly equivalent to $D_i$ on $X_i$. Then we have 
$\rho_{*}(D)  = \underline{d}$ as in Case 1. 

For $i = 1, 2$, 
we set 
\[
  \varepsilon_i = 
  \begin{cases} 
   1 & \quad \textup{(if $p \in {\rm Bs}(D_i)$)},  \\
   0 & \quad  \textup{(if $p \not\in {\rm Bs}(D_i)$)}, 
   \end{cases}
\]
so that $r_{X_i}(D_i^\prime) = r_{X}(D_i) - (1- \varepsilon_i)$. 
Then it follows from Lemma~\ref{lemma:Inaba} that 
\begin{align*}
  r_X(D) 
  & =  r_{X_1}(\rest{D}{X_1}) + r_{X_2}(\rest{D}{X_2}) 
  + \varepsilon_1 \varepsilon_2 \\
  & = r_{X_1}(D_1) + r_{X_2}(D_2) + \varepsilon_1 \varepsilon_2 \\
  & =  r_{X_1}(D_1^\prime) + r_{X_2}(D_2^\prime) + (1- \varepsilon_1)
  + (1- \varepsilon_2) +   \varepsilon_1 \varepsilon_2
   \\ 
  & \geq r_{\bar{G}_1}(\underline{d_1} - [v_1]) + 
 r_{\bar{G}_2}(\underline{d_2} - [v_2])  + 1 + (1- \varepsilon_1)(1- \varepsilon_2) \\
 & \geq r_{\bar{G}_1}(\underline{d_1} - [v_1]) + 
 r_{\bar{G}_2}(\underline{d_2} - [v_2])  + 1. 
\end{align*}
This gives the desired properties in the remaining case, 
thus completing the proof. 
\QED 

Next, we reinterpret Proposition~\ref{prop:main}. 

\begin{Lemma}
\label{lemma:reinterpret}
In the setting of Proposition~\ref{prop:main}, 
we have
%\[ 
\begin{equation}
\label{eqn:claim}
 r_{\bar{G}}(\underline{d}) \leq 
 \begin{cases}
 r_{\bar{G}_1}(\underline{d_1} - [v_1]) + 
 r_{\bar{G}_2}(\underline{d_2} - [v_2])  + 2 &  
 \text{\textup{(}if $v_i \not\in {\rm Bs}(\left|\underline{d_i}\right|^\bullet)$ for each $i = 1, 2$\textup{)}, } \\
 r_{\bar{G}_1}(\underline{d_1} - [v_1]) + 
 r_{\bar{G}_2}(\underline{d_2} - [v_2]) + 1 & 
 \text{\textup{(}otherwise\textup{)}}. 
 \end{cases}
\end{equation}
%\]
\end{Lemma}

\Proof
First we consider the case where $v_i \not\in {\rm Bs}(\left|\underline{d_i}\right|^\bullet)$ for each $i = 1, 2$. 
%This means that $\left|\underline{d_i}\right|^\bullet \neq \emptyset$ and 
Then, since
$r_{\bar{G}_i}(\underline{d_i} - [v_i]) 
= r_{\bar{G}_i}(\underline{d_i}) -1$,
it follows from Proposition~\ref{prop:main} that 
\[
 r_{\bar{G}}(\underline{d}) 
 \leq 
 r_{\bar{G}_1}(\underline{d_1}) + 
 r_{\bar{G}_2}(\underline{d_2}) 
 =  r_{\bar{G}_1}(\underline{d_1} - [v_1]) + 
 r_{\bar{G}_2}(\underline{d_2} - [v_2])  + 2. 
\]
This shows the inequality (\ref{eqn:claim}) in this case.

Next consider 
the case where $v_1\in {\rm Bs}(\left|\underline{d_1}\right|^\bullet)$ 
and $v_2 \not\in {\rm Bs}(\left|\underline{d_2}\right|^\bullet)$. (The case of 
$v_1 \not\in {\rm Bs}(\left|\underline{d_1}\right|^\bullet)$ 
and $v_2 \in {\rm Bs}(\left|\underline{d_2}\right|^\bullet)$ 
is shown in the same way.) Then 
we have $r_{\bar{G}_1}(\underline{d_1} - [v_1]) 
= r_{\bar{G}_1}(\underline{d_1})$, 
%$\left|\underline{d_2}\right|^\bullet \neq \emptyset$ 
and $r_{\bar{G}_2}(\underline{d_2} - [v_2]) 
= r_{\bar{G}_2}(\underline{d_2}) -1$. It follows from 
Proposition~\ref{prop:main} that 
\[
 r_{\bar{G}}(\underline{d}) 
 \leq 
 r_{\bar{G}_1}(\underline{d_1}) + 
 r_{\bar{G}_2}(\underline{d_2}) 
 =  r_{\bar{G}_1}(\underline{d_1} - [v_1]) + 
 r_{\bar{G}_2}(\underline{d_2} - [v_2])  + 1,
\]
which shows (\ref{eqn:claim}) in this case.

Finally consider the case where $v_i\in {\rm Bs}(\left|\underline{d_i}\right|^\bullet)$ 
for each $i = 1, 2$. 
This means that $r_{\bar{G}_i}(\underline{d_i} - [v_i]) 
= r_{\bar{G}_i}(\underline{d_i})$, and 
Proposition~\ref{prop:main} gives 
\[
 r_{\bar{G}}(\underline{d}) 
 \leq 
 r_{\bar{G}_1}(\underline{d_1}) + 
 r_{\bar{G}_2}(\underline{d_2}) + 1
 =  r_{\bar{G}_1}(\underline{d_1} - [v_1]) + 
 r_{\bar{G}_2}(\underline{d_2} - [v_2])  + 1. 
\]
Thus we obtain \eqref{eqn:claim}. 
\QED

To prove Theorem~\ref{thm:cor:main}, 
we consider the 
following condition for a vertex-weighted graph~$\bar{G}$. 

\bigskip
\begin{enumerate}
\item[(FS)]
There exists a nodal curve $X$ with dual graph $\bar{G}$ such that, for any 
$\underline{d} \in \Div(G)$, there exists a Cartier divisor $D$ on 
$X$ such that $\rho_{*}(D) = \underline{d}$ 
and $r_X(D) \geq r_{\bar{G}}(\underline{d})$, where 
$\rho_{*}: \Div(X) \to \Div(G)$ is the map defined in \eqref{eqn:sp}. 
\end{enumerate}

\bigskip
We note that if a vertex-weighted graph $\bar{G}$ satisfies the condition 
(FS), then we have $r_{\bar{G}}^{\alg}(\delta) \geq r_{\bar{G}}(\delta)$ 
for any divisor class $\delta \in \Pic(G)$. Indeed, let $\delta$ be any divisor class of 
$G$. 
We take the nodal curve $X$ in the condition (FS). Then for any representative $\underline{d} \in \Div(G)$ 
of  $\delta$, we take a Cartier 
divisor $D$ on $X$ as in (FS). 
With $X$ as above, we obtain 
$
  \min_{\underline{d} \in \delta} \left\{
  \max_{L \in \Pic^{\underline{d}}(X)} \{r_X(L)\}
  \right\} \geq r_{\bar{G}}(\delta)
$. Thus we get 
$r_{\bar{G}}^{\alg}(\delta) \geq r_{\bar{G}}(\delta)$. 

\medskip
{\bf Proof of Theorem~\ref{thm:cor:main}.}\quad
Let $\bar{G}$ be a hyperelliptic vertex-weighted graph. 
We will show that $\bar{G}$ satisfies the condition (FS) 
by the induction on the number of bridges. As is explained 
as above, we will then have the desired inequality 
$r_{\bar{G}}^{\alg}(\delta) \geq r_{\bar{G}}(\delta)$ for any divisor class 
$\delta$ on $G$.

If $G$ has no bridges, then \cite[Proposition~1.5 and its proof and Theorem~8.2]{KY} tells us that $\bar{G}$ satisfies the condition (FS) 
(More generally, if there are at most $(2 \omega(v)+2)$ ``positive-type'' bridges emanating from each vertex $v\in V(G)$, then $\bar{G}$ satisfies the condition (FS): See \cite{KY}.) Also, 
a vertex-weighted graph of genus at most $1$ satisfies the condition (FS) (cf. \cite[Proposition~1.5 and its proof and Proposition~7.5]{KY}). 

\smallskip
Now we consider the general case, and suppose that $\bar{G}$ has a bridge. Let $G_1$ and $G_2$ be the connected components
of $G\setminus\{e\}$, and 
set $\bar{G}_i = (G_i, \rest{\omega}{V(G_i)})$ for $i = 1, 2$. 
Then we find that $\bar{G}_i$ is a hyperelliptic or $g(\bar{G}_i) \leq 1$ (see \cite[\S5.2]{BN2} or \cite[Lemma~3.4]{KY}). 

By the induction on the the number of bridges, 
we may and do assume that $\bar{G}_i$ 
satisfies the condition (FS) for each $i = 1, 2$. 
Thus there exists 
a nodal curve $X_i$  such that, for any 
$\underline{e_i} \in \Div(G_i)$, there exists a Cartier 
divisor $E_i$ on $X_i$ satisfying $\rho_{i*}(E_i) = \underline{e_i}$ 
and $r_{X_i}(E_i) \geq r_{\bar{G}_i}(\underline{e_i})$, 
where $\rho_{i*}: \Div(X_i) \to \Div(G_i)$ is the map defined in \eqref{eqn:sp}. 

Let $p_i$ be a smooth point of $X_i$ for
each $i = 1,2$.
Then we patch $X_1$ and $X_2$ by 
$p_1 = p_2 \, (=: p)$ to obtain a nodal curve $X$
such that $X = X_1 \cup X_2$ and $X_1 \cap X_2 = \{ p \}$.
Here we take
each $p_i$
so that
$X = X_1 \cup X_2$ is a nodal curve 
with dual graph $\bar{G}$ 
and that each $G_i$ is the subgraph of $G$
corresponding to the component $X_i$.
Let $\rho_{*}: \Div(X) \to \Div(G)$ be the map defined in \eqref{eqn:sp}. 

\smallskip
We prove that, with this $X$, $\bar{G}$ satisfies the condition (FS). Indeed, let $\underline{d}$ be any divisor on $G$. 
For $i = 1, 2$, let 
$\underline{d}_i$ be the restriction of $\underline{d}$ to $G_i$.  
By Lemma~\ref{lemma:useful}, there exists a Cartier divisor $D$ on 
$X$ satisfying $\rho_{*}(D) = \underline{d}$ and 
\[
 r_X(D) \geq 
 \begin{cases}
 r_{\bar{G}_1}(\underline{d_1} - [v_1]) + 
 r_{\bar{G}_2}(\underline{d_2} - [v_2])  + 2 &  
 \text{\textup{(}if $v_i \not\in {\rm Bs}(\left|\underline{d_i}\right|^\bullet)$ for each $i = 1, 2$\textup{)}, } \\
 r_{\bar{G}_1}(\underline{d_1} - [v_1]) + 
 r_{\bar{G}_2}(\underline{d_2} - [v_2]) + 1 & 
 \text{\textup{(}otherwise\textup{)}}. 
 \end{cases}
\]
By Lemma~\ref{lemma:reinterpret}, the right-hand side is at least $r_{\bar{G}}(\underline{d})$. 
Thus we obtain $r_X(D)  \geq r_{\bar{G}}(\underline{d})$, which shows that 
$\bar{G}$ satisfies the condition (FS). 
\QED

%%%%%%%%%%%%%%%%%%%%%%%%%%%%%%
% Rank of divisors on graphs and curves of genus $3$ %
%%%%%%%%%%%%%%%%%%%%%%%%%%%%%%
\setcounter{equation}{0}
\section{Rank of divisors on graphs and curves of genus $3$}
\label{sec:genus:3}
In this section, we prove Proposition~\ref{prop:KY:3} and 
then Theorem~\ref{thm:main:2}. 
Let $R$ be a complete discrete valuation ring with fractional field $\KK$ and 
residue field $k$. Let $\Xscr$ be a regular, generically smooth, semi-stable $R$-curve. Let 
$\bar{G} = (G, \omega)$ be the reduction graph of $\Xscr$. 
Let $\widetilde{\rho}_*: \Div(\Xscr_\KK) \to \Div(G)$ be 
the specialization map defined in \eqref{eqn:sp:2}.

We begin the proof of Proposition~\ref{prop:KY:3}. 

\medskip
{\bf Proof of Proposition~\ref{prop:KY:3}}.\quad 
Recall that $\bar{G} = (G, \omega)$ is a non-hyperelliptic graph of genus $3$ and that $\Xscr$ is a regular, generically smooth, semi-stable $R$-curve $\Xscr$ 
with reduction graph $\bar{G}$. 

First we claim that, if $\deg(\underline{d}) \leq 2$, then 
$r_{\bar{G}}(\underline{d}) \leq 0$. Indeed, to argue by contradiction, suppose that 
$r_{\bar{G}}(\underline{d}) \geq 1$. Since $r_{\bar{G}}(\underline{d}) \leq \deg (\underline{d})$, this means (A) $\deg (\underline{d}) = 1$ and 
$r_{\bar{G}}(\underline{d}) = 1$; (B) $\deg (\underline{d}) = 2$ and 
$r_{\bar{G}}(\underline{d}) = 1$; or (C) $\deg (\underline{d}) = 2$ and 
$r_{\bar{G}}(\underline{d}) = 2$. In (A),  the existence of $\underline{d}$ forces $\bar{G}^\bullet$ to be a tree, 
which is a contradiction. In (B), 
$\bar{G}$ is hyperelliptic, which is excluded at the beginning. 
In (C), there exists a vertex $v$ of $\bar{G}^\bullet$
such that $r_{\bar{G}^\bullet}(\underline{d} - [v]) = 1$ and 
$\deg(\underline{d} - [v]) = 1$, where we regard 
$\underline{d} \in \Div(\bar{G}^\bullet)$ via the natural 
embedding 
$\Div(G) \subseteq \Div(\bar{G}^\bullet)$ 
as before. 
The existence of the divisor $\underline{d} - [v]$ 
forces $\bar{G}^\bullet$ to be a tree, 
which is a contradiction. 
Hence we obtain the claim. 

\medskip
{\bf Case 1.}\quad 
Suppose that $r_{\bar{G}}(\underline{d}) = -1$.
By the surjectivity of the homomorphism \eqref{eqn:surj:rho},
there exists $\widetilde{D} \in \Div(\Xscr_{\KK})$ with 
$\widetilde{\rho}_*(\widetilde{D}) = \underline{d}$.
% (cf. \eqref{eqn:surj:rho}). 
Then the specialization lemma (Theorem~\ref{thm:AC}) tells us 
that $-1 = r_{\bar{G}}(\underline{d}) 
\geq r_{\Xscr_{\KK}}(\widetilde{D})$. It follows that 
$r_{\Xscr_{\KK}}(\widetilde{D}) = -1$.
Thus
the equality 
$r_{\Xscr_\KK}(\widetilde{D}) = r_{\bar{G}}(\underline{d})$ holds. 

\medskip
{\bf Case 2.}\quad 
Suppose that $r_{\bar{G}}(\underline{d}) = 0$.
Then it follows from Lemma~\ref{lemma:for:genus:3} that 
there exists an effective divisor $\underline{e} \in \Div(G)$ such that 
$\underline{e}$ is linearly equivalent to $\underline{d}$ in $G$. 
Since the homomorphism \eqref{eqn:surj:rho}
induces a surjective map between the sets of effective divisors,
there exists an effective divisor $\widetilde{E} \in \Div(\Xscr_{\KK})$ 
with $\widetilde{\rho}_*(\widetilde{E}) = \underline{e}$.
% (cf. \eqref{eqn:surj:rho}).  
Now we use Raynaud's theorem (Theorem~\ref{thm:Raynaud} below) as in 
the proof of \cite[Theorem~1.5]{KY}. 
It follows that there exists a principal divisor $\widetilde{N} \in \Div(\Xscr_{\KK})$
such that $\widetilde{\rho}_*(\widetilde{N}) = \underline{d} - \underline{e}$. 
We set  $\widetilde{D} =  \widetilde{E} + \widetilde{N}$. Then 
$\rho_*(\widetilde{D}) = \underline{d}$. 
Since $\widetilde{D}$ is linearly equivalent to $\widetilde{E}$ on $\Xscr_\KK$, 
we have 
$
r_{\Xscr_\KK}(\widetilde{D}) 
= r_{\Xscr_\KK}(\widetilde{E}) \geq 0. 
$
Since 
\[
0 = r_{\bar{G}} (\underline{d}) = r_{\bar{G}}(\underline{e}) 
\geq r_{\Xscr_{\KK}}(\widetilde{E})
\]
by the specialization lemma (Theorem~\ref{thm:AC}), we obtain the equality 
$r_{\Xscr_\KK}(\widetilde{D}) = r_{\bar{G}}(\underline{d})\; (= 0)$. 

\medskip
{\bf Case 3.}\quad 
Suppose that $r_{\bar{G}}(\underline{d}) \geq 1$. 
By the above claim, we have $\deg (\underline{d}) \geq 3$. 
We put $\underline{d}^\prime := K_{\bar{G}} - \underline{d} \in \Div(G)$. 
Then $\deg (\underline{d}^\prime) = 4 - \deg (\underline{d}) \leq 1$. 
Thus $r_{\bar{G}} (\underline{d}^\prime) \leq \deg (\underline{d}^\prime) \leq 1$. Since $\bar{G}^\bullet$ is not a tree, we have $r_{\bar{G}} (\underline{d}^\prime) \neq 1$. It follows that $r_{\bar{G}} (\underline{d}^\prime) \leq 0$.  
By Cases 1 and 2, 
there exists a divisor $\widetilde{D'} \in \Div (\Xscr_\KK)$ 
such that $\widetilde{\rho}_*(\widetilde{D'}) = \underline{d}^\prime$ and 
$r_{\Xscr_\KK}(\widetilde{D'}) = r_{\bar{G}} (\underline{d}^\prime) $.

By \cite[Remark~4.18~and~Remark~4.21]{B}, there exists a canonical divisor $K_{\Xscr_\KK}$ of $\Xscr_\KK$ such that $\widetilde{\rho}_*(K_{\Xscr_\KK}) = K_{\bar{G}}$. 
We set $\widetilde{D} := K_{\Xscr_\KK} - \widetilde{D'}$.  
Then we have 
$\widetilde{\rho}_*(\widetilde{D}) = K_{\bar{G}} - \underline{d}^\prime =\underline{d}$.
Further, the Riemann--Roch formulae on $\Xscr_\KK$ and 
$\bar{G}$ (cf. \cite[Theorem~3.8]{AC}) give 
\[
r_{\Xscr_\KK}(\widetilde{D})
= - 2 + \deg (\widetilde{D}) + r_{\Xscr_\KK}(\widetilde{D'})
= - 2 + \deg (\underline{d}) + r_{\bar{G}}(\underline{d}')
=
r_{\bar{G}}(\underline{d})
.
\]
Thus we obtain Proposition~\ref{prop:KY:3}. 
\QED

\medskip
{\bf Proof of Theorem~\ref{thm:main:2}.}\quad
The proof goes in the same way as in \cite{KY}; 
Theorem~\ref{thm:main:2} will be deduced from Proposition~\ref{prop:KY:3}.
%by the same argument in  \cite{KY} (that is due to Caporaso). 

Recall that $k$ is a fixed algebraically closed field. We take 
a complete discrete valuation ring $R$ with residue field $k$. For example, 
we may take $R$ as the ring of formal power series $k [\![ t ]\!]$ over $k$. 
Let $\KK$ be the fractional field of $R$. 
We take a regular, generically smooth, semi-stable $R$-curve $\Xscr$ 
with reduction graph $\bar{G}$. We note that such $\Xscr$ always exists:  
See \cite[Theorem~B.2]{B}.  

Let $\Xscr_{\KK}$ denote the generic fiber of $\Xscr$, and 
$X$ the special fiber of 
$\Xscr$. For $\underline{d} \in \Div(G)$, Proposition~\ref{prop:KY:3} shows that there exists a divisor $\widetilde{D} 
\in \Div(\Xscr_{\KK})$ such that 
$\widetilde\rho_*(\widetilde{D}) = \underline{d}$ and 
$r_{\bar{G}}(\underline{d}) = r_{\Xscr_\KK}(\widetilde{D})$. 
Let $\widetilde{\Dscr}$ be the Zariski closure of $\widetilde{D}$ in 
$\Xscr$. We denote by $\OO_{\Xscr}(\widetilde{\Dscr})$ the 
invertible sheaf on $\Xscr$ associated to $\widetilde{\Dscr}$. 
Let $D \in \Div(X)$ be a divisor on $X$ such that 
the associated invertible sheaf $\OO_X(D)$ is isomorphic to 
$\rest{\OO_{\Xscr}(\widetilde{\Dscr})}{X}$. By the upper-semicontinuity 
of the cohomology, we have $r_{X}(D) \geq r_{\Xscr_\KK}(\widetilde{D})$. 
Hence $r_{X}(D) \geq r_{\bar{G}}(\underline{d})$. Also, by \eqref{eqn:compatibility}, we have $\rho_*(D) = \underline{d}$. 
If follows that $\bar{G}$ satisfies the condition (FS) in Section~\ref{sec:bridge}, and we get 
$r_{\bar{G}}^{\alg}(\delta) \geq r_{\bar{G}}(\delta)$ for any 
divisor class $\delta \in \Pic(G)$. 
\QED
 
\bigskip 
In the rest of this section, 
we will show a metric graph version of Proposition~\ref{prop:KY:3}. 
Let $\Xscr$ be a regular, generically smooth, semi-stable $R$-curve with reduction graph 
$\bar{G} = (G, \omega)$. Let $\Gamma$ be the metric graph 
associated to $G$, 
where each edge of $G$ is assigned length one. 
Let $\Gamma_\QQ$ be the the set of points of $\Gamma$ whose distance 
from every vertex of $G$ is rational. 

We follow the arguments in \cite[Section~2.3]{B}. 
Let $\KK^\prime/\KK$ be a finite extension. Let $R^\prime$ be 
the ring of integers of $\KK^\prime$. 
Then $R^\prime$ is a complete discrete valuation ring with residue field $k$. 
Let $\Xscr^\prime$ be the minimal resolution of 
$\Xscr\times_{\Spec(R)} \Spec(R^\prime)$. Then 
$\Xscr^\prime$ is a regular, generically smooth, semi-stable $R^\prime$-curve 
with generic fiber $\Xscr\times_{\Spec(\KK)} \Spec(\KK^\prime)$.  
Let $e(\KK^\prime/\KK)$ be the ramification index of  $\KK^\prime/\KK$. 
The dual graph $\bar{G^\prime} = (G^\prime, \omega^\prime)$ 
of the special fiber of $\Xscr^\prime$ is 
the graph obtained by inserting $e(\KK^\prime/\KK)-1$ vertices 
to each edge of $G$, and $\omega^\prime$ is the extension of $\omega$, where $\omega^\prime(w) = 0$ for any $w \in V(G^\prime)\setminus V(G)$. 
If we assign a length of $1/e(\KK^\prime/\KK)$ to each edge of $G^\prime$, 
then the corresponding metric graph equals $\Gamma$. 
The pair $\Gamma$ with a vertex-weight function $\Gamma\to\ZZ$ (given by the zero extension of $\omega$) is denoted by $\bar{\Gamma}$

Let $\bar{\KK}$ be an algebraic closure of $\KK$. 
For $\widetilde{D} \in \Div(\Xscr_{\bar{\KK}})$, we take a finite extension 
$\KK^\prime/\KK$ such that $\widetilde{D} \in \Div(\Xscr(\KK^\prime))$, 
and then we set $\widetilde{\tau}_*(\widetilde{D}) = \widetilde{\rho^\prime}_*(\widetilde{D})$, 
where $\widetilde{\rho^\prime}_*$ is the specialization map for $\Xscr^\prime$. 
This gives rise to the specialization map 
\[
  \widetilde{\tau}_*: \Div(\Xscr_{\bar{\KK}}) \to \Div(\Gamma_\QQ). 
\]
(This map is denoted by $\tau_*$
in \cite{KY}.
Here we write $\widetilde{\tau}_*$ instead
because of the compatibility 
with the notation $\widetilde{\rho}_*$; cf. Remark~\ref{rmk:notation}.)

For each $\underline{d} \in \Div(\Gamma_\QQ)$, we take a graph 
$\bar{G}^\prime = (G^\prime, \omega^\prime)$ with $\Supp(\underline{d}) \subset V(G^\prime)$, and define $r_{\bar{\Gamma}}(\underline{d}) 
:= r_{\bar{G^\prime}}(\underline{d})$, 
which does not depend on the choice of 
$\bar{G}^\prime$ by \cite[\S1]{AC}. 
By Amini--Caporaso's specialization lemma (Theorem~\ref{thm:AC}), 
we have 
$r_{\bar{G^\prime}}(\rho^\prime_*(\widetilde{D})) \geq r_{\Xscr_{\KK^\prime}}(\widetilde{D})$ for any $\widetilde{D} \in \Div(\Xscr_{\KK^\prime})$. 
As is mentioned in \cite[Remark~2.9]{B}, if 
$\widetilde{D} \in \Div(\Xscr_{\KK^\prime}) \setminus \Div(\Xscr(\KK^\prime))$, then $\widetilde{\rho^\prime}_*(\widetilde{D})$ and $\widetilde{\tau}_*(\widetilde{D})$ may be different, 
but $\widetilde{\rho}_*(\widetilde{D})$ 
and $\widetilde{\tau}_*(\widetilde{D})$ are at least linearly equivalent in $G^\prime$. 
Thus, we have the specialization lemma for vertex-weighted metric graph: 
For any $\widetilde{D} \in \Div(\Xscr_{\bar{\KK}})$, one has 
\[
  r_{\bar{\Gamma}}(\tau_*(\widetilde{D})) \geq r_{\Xscr_{\bar{\KK}}}(\widetilde{D}). 
\]

Also for metric graphs, we have Raynaud's theorem,
which asserts the 
surjectivity of the map $\rest{\widetilde{\tau}_*}{\Prin(\Xscr_{\bar{\KK}})}: 
\Prin(\Xscr_{\bar{\KK}}) \to \Prin(\Gamma_\QQ)$ (see \cite[Corollary~A.9]{B}). 
Then, by the same argument as in the proof of Proposition~\ref{prop:KY:3}, 
we obtain the following proposition. 

\begin{Proposition}
\label{prop:KY:3:2}
Let $R$ be a complete discrete valuation ring 
with fractional field $\KK$ and residue field $k$. 
Let $\bar{G} = (G , \omega)$ be a non-hyperelliptic graph of 
genus $3$, and $\Gamma$ the metric graph associated to $G$, 
where each edge of $G$ is assigned length one. Let $\Xscr$ be a regular, generically smooth, semi-stable $R$-curve with 
reduction graph $\bar{G}$.  Then the following condition 
\textup{(C)}  
holds. 
\begin{enumerate}
\item[(C)]
For any $\underline{d} \in \Div (\Gamma_\QQ)$, there exists a divisor $\widetilde{D} \in \Div (\Xscr_{\bar{\KK}})$
such that $\widetilde{\tau}_*(\widetilde{D}) = \underline{d}$
and $r_{\bar{\Gamma}}(\underline{d}) = r_{\Xscr_{\bar{\KK}}}(\widetilde{D})$.
\end{enumerate}
\end{Proposition}

\begin{Remark}
\label{rmk:last}
Let $R$ be a complete discrete valuation ring 
with fractional field $\KK$ and residue field $k$. 
Let $\bar{G} = (G, \omega)$ be a vertex-weighted graph. In \cite{KY}, we have asked  
under what condition on $\bar{G}$ there exists a regular, generically smooth, semi-stable $R$-curve with reduction graph $\bar{G}$ that satisfies the conditions (F) and (C) in Proposition~\ref{prop:KY:3} and Proposition~\ref{prop:KY:3:2}. In \cite{KY}, when $\ch(k) \neq 2$, we have completely answered this question for hyperelliptic graphs: A hyperelliptic graph $\bar{G} =  (G, \omega)$ satisfies the conditions (F) and (C) if and only if 
every vertex $v$ of $G$ has at most $2 \omega(v) + 2$ ``positive-type'' bridges 
emanating from it. In this paper, we answer this question for non-hyperelliptic graphs of genus $3$:  Every non-hyperelliptic graph of genus $3$ satisfies (F) and (C). It is then natural to ask this question for non-hyperelliptic graphs of genus $4$. In this case, the arguments in the proof 
of Proposition~\ref{prop:KY:3} show the existence of a desired lift $\widetilde{D}$ of $\underline{d}$ except for 
divisors $\underline{d}$ with 
$\deg(\underline{d}) = 3$ and $r_{\bar{G}}(\underline{d})= 1$.   
\end{Remark}

\begin{Remark}
\label{rmk:last:2}
Assume that $\ch(k) \neq 2$. 
Let $\bar{G} = (G, \omega)$ be a vertex-weighted graph. Let 
$e_1, \ldots, e_r$ 
be the set of bridges of $G$, and we write 
$G \setminus \{e_1, \ldots, e_r\} = G_1 \cup \cdots \cup G_{r+1}$ as the disjoint union of connected finite graphs. We set 
$\bar{G}_i = (G_i, \rest{\omega}{G_i})$. The proofs of Theorem~\ref{thm:cor:main} and 
Theorem~\ref{thm:main:2} show that 
hyperelliptic graphs and graphs of genus at most $3$ satisfy 
the condition (FS). It follows from the proof of Theorem~\ref{thm:cor:main} that if each $\bar{G}_i$ is hyperelliptic or 
of genus at most $3$, then  $\bar{G}$ satisfies 
the condition (FS), and thus we have 
$r_{\bar{G}}^{\alg}(\delta) \geq r_{\bar{G}}(\delta)$ 
for any divisor class $\delta \in \Pic(G)$.
\end{Remark}

%%%%%%%
% Appendix %
%%%%%%%
\renewcommand{\theTheorem}{A.\arabic{Theorem}}
\renewcommand{\theClaim}{A.\arabic{Theorem}.\arabic{Claim}}
\renewcommand{\theequation}{A.\arabic{equation}}
\renewcommand{\thesubsection}{A.\arabic{subsection}}
\setcounter{Theorem}{0}
\setcounter{subsection}{0}
\setcounter{Claim}{0}
\setcounter{equation}{0}
\section*{Appendix: Raynaud's theorem}
The purpose of this appendix is to show that, for a finite graph with loops, 
the specialization map between principal divisors is still surjective. 
Our proof of the surjectivity will be given by reducing to the case of 
loopless finite graphs. 
The surjectivity in the loopless case
is shown in Baker \cite{B}. 

In \cite{B},
the surjectivity of the specialization map
(in the loopless case)
is attributed to Raynaud
because this surjectivity follows 
from re-interpretation of Raynaud's results in \cite{Ra} (see \cite[Appendix~A]{B}). 
In this paper, we also call
Theorem~\ref{thm:Raynaud},  which asserts the surjectivity, 
Raynaud's theorem. 

Let $k$ be an algebraically closed field as before. Let $R$ be a complete valuation ring with residue field $k$. Let $\KK$ be the fractional field of $R$. 
Let $\Xscr \to \Spec (R)$ be a 
regular, generically smooth, semi-stable $R$-curve.
We write $X$ for the special fiber of $\Xscr$, and 
$\bar{G} = (G, \omega)$ for the dual graph of $X$. 
Let $\widetilde{\rho}_*: \Div (\Xscr_\KK) \to \Div (G)$ be the specialization map  
defined in \eqref{eqn:sp:2}.

\begin{Theorem} 
\label{thm:Raynaud}
The specialization map between principal divisors is surjective. Namely, 
$\rest{\widetilde{\rho}_*}{\mathrm{Prin} (\Xscr_\KK)} :\mathrm{Prin} (\Xscr_\KK) \to \mathrm{Prin} (G)$
is surjective.
\end{Theorem}

\Proof
We put $p := \ch (k) \geq 0$. 
When $G$ is loopless, then the assertion is exactly 
\cite[Corollary~A.8]{B}. We will reduce the general case to the loopless case. 

Let $d$ be an integer with
$d \geq 2$.
When $p >0$, we 
require that $(d , p) = 1$.
We fix a finite Galois extension $\KK^\prime$ 
of $\KK$ of degree $d$. (For example, we may take 
$\KK^\prime = \KK(\sqrt[d]{\pi})$, where $\pi\in R$ is 
a uniformizer of $R$.)  
Since $k$ is algebraically closed and $\KK^\prime/\KK$ is 
a Galois extension of degree $d$, 
the ramification index $e(\KK^\prime/\KK)$ equals $d$. 
We denote by $R^\prime$ the ring of integers of $\KK^\prime$.

Let $\Xscr^\prime$ be the minimal resolution of 
$\Xscr\times_{\Spec(R)}\Spec(R^\prime)$. 
Let $\nu: \Xscr^\prime \to \Xscr$ be the natural map. By slight abuse of notation, we denote the restriction of $\nu$ to the generic fibers by the same notation $\nu$. Let $X^\prime$ be the special fiber of $\Xscr^\prime$. 
Let $G^\prime$ be dual graph of 
$X^\prime$,
and let $\widetilde{\rho}^\prime_* : \Div (\Xscr_{\KK^\prime}) \to \Div (G^\prime)$
be the specialization map with respect to $\Xscr^\prime$. 
Since $G^\prime$ is the graph obtained by
inserting $(d-1)$ vertices to each edge of $G$,
we have a natural embedding $V(G) \subseteq V(G^\prime)$ and also  
$\Div(X) \subseteq \Div(X^\prime)$.  

\begin{Claim} \label{claim:surj-prins-2}
Let $\widetilde{D} \in \Div (\Xscr_\KK)$ such that  
$\widetilde{\rho}_*^\prime  (\nu^* (\widetilde{D})) \in \Div (G)$. Then
$\widetilde{\rho}_*^\prime (\nu^* (\widetilde{D})) = \widetilde{\rho}_*  (\widetilde{D})$.
\end{Claim}

Indeed, we take any $v \in V(G)$.
Let $C_v$ be the irreducible component of 
$X$ corresponding to $v$
and let $C_v^\prime$ be the irreducible component of $X^\prime$
with $\nu (C_v^\prime) = C_v$.
Let $\mathscr{D} \in \Div ( \Xscr )$ 
and $\mathscr{D}^\prime \in \Div ( \Xscr^\prime )$ be 
the Zariski closures of $\widetilde{D}$ and $\nu^*(\widetilde{D})$ respectively.
We have $\widetilde{\rho}_*(\widetilde{D}) (v) = (C_v \cdot \mathscr{D})$
and $\widetilde{\rho}_*^\prime(\nu^* (\widetilde{D})) (v) = (C^\prime_v \cdot \mathscr{D}^\prime)$.

Since $\nu_* (\mathscr{D}^\prime) = d \mathscr{D}$,
we have
$d (C_v \cdot \mathscr{D})
= (\nu^* (C_v) \cdot \mathscr{D}^\prime)$ by the projection formula.
By the assumption that
$\widetilde{\rho}_*^\prime (\nu^* (D)) \in \Div (G)$,
we have $(E^\prime \cdot \mathscr{D}^\prime) = 0$ for any exceptional prime divisor $E^\prime$
for $\nu$.
Since $\nu^* (C_v) - d C_v^\prime$ is a linear combination of exceptional divisors,
it follows that 
$(\nu^* (C_v) \cdot \mathscr{D}^\prime) = d( C_v^\prime \cdot \mathscr{D}^\prime)$.

Then 
\[
\widetilde{\rho}_*(\widetilde{D}) (v) = (C_v \cdot \mathscr{D})
= \frac{(\nu^* (C_v) \cdot \mathscr{D}^\prime)}{d}
= (C_v^\prime \cdot \mathscr{D}^\prime) = \widetilde{\rho}^\prime( \nu^{*} (\widetilde{D})) (v)
.
\]
Since $\widetilde{\rho}_*^\prime  (\nu^* (\widetilde{D})) \in \Div (G)$ and 
$v \in V(G)$ is arbitrary, we obtain Claim~\ref{claim:surj-prins-2}. 

\medskip
Let $\sigma_{1} , \ldots , \sigma_{d}$ 
be the elements of $\Gal(\KK^\prime / \KK)$.
Each $\sigma_{i}$ induces an automorphism $\sigma_{i}^* : \Spec (R^\prime) \to \Spec (R^\prime)$,
and
an automorphism
$\varphi_{i} : \Xscr^\prime \to \Xscr^\prime$
over $R$ (induced from the cartesian product). 
Let $\varphi_{i}^* : \Div (\Xscr^\prime) \to \Div (\Xscr^\prime)$ and 
$\varphi_{i}^* : \Div (\Xscr_{\KK^\prime}) \to \Div (\Xscr_{\KK^\prime})$ be the induced maps. 

\begin{Claim} 
\label{claim:surj-prins-1}
For any $\widetilde{D}^\prime \in \Div(\Xscr_{\KK^\prime})$, we have 
$\widetilde{\rho}_*^\prime\left((\varphi_{i})^*(\widetilde{D}^\prime)\right) = \widetilde{\rho}_*^\prime(\widetilde{D}^\prime)$ for $i = 1, \ldots, d$. 
\end{Claim}

Indeed, since $\sigma_i$ induces the trivial action on the residue field $k$,
the restriction of $\varphi_i$ to the special fiber $X^\prime$  
is trivial.
Thus $(\varphi_i)_* (C^\prime) = C^\prime$
for any irreducible component $C^\prime$ of $X^\prime$. 

We take any $\widetilde{D}^\prime \in \Div (\Xscr_{\KK^\prime})$ and let $\mathscr{D}^\prime$ be
the Zariski closure of $\widetilde{D}^\prime$ in $\Xscr^\prime$.
Note that $\varphi_i^* (\mathscr{D}^\prime)$ is
the Zariski closure of $\varphi_i^* (\widetilde{D}^\prime)$.
For any $v \in V (G^\prime)$, 
let $C^\prime_v$ be the corresponding irreducible component
of $X^\prime$.
Then 
\[
\widetilde{\rho}_*^\prime (\varphi_i^* (\widetilde{D}^\prime)) (v)
=
(C^\prime_v \cdot \varphi_i^*(\mathscr{D}^\prime)) =
((\varphi_i)_* ( C^\prime_v) \cdot \mathscr{D}^\prime) = (C^\prime_v \cdot \mathscr{D}^\prime)
=
\widetilde{\rho}_*^\prime (\widetilde{D}^\prime) (v)
,
\]
which shows the desired equality.
We obtain Claim~\ref{claim:surj-prins-1}. 

\medskip
We take any $\underline{n} \in \Prin(G)$.
Then
$\underline{n} \in \mathrm{Prin} (G^\prime)$.
Since $G^\prime$ is loopless, we know that 
$\widetilde{\rho}_*^\prime : \mathrm{Prin} (\Xscr_{\KK^\prime}) \to \mathrm{Prin} (G^\prime)$ is 
surjective
by \cite[Corollary~A.8]{B}. 
Let $f$ be a non-zero rational function on $\Xscr_{\KK^\prime}$
such that $\widetilde{\rho}_*^\prime (\zero(f)) = \underline{n}$.
We set $g^\prime := \varphi_1^* (f) \cdots \varphi_d^* (f)$,
which is a non-zero rational function on $\Xscr_{\KK^\prime}$.
Then $\zero (g^\prime) = \varphi_1^* (\zero ( f)) + 
\cdots + \varphi_d^* (\zero ( f ))$,
so that
Claim~\ref{claim:surj-prins-1}
tells us 
that $\widetilde{\rho}_*^\prime ( \zero (g^\prime) ) = d\, \underline{n}$.
Since $g^\prime$ is a $\Gal(\KK^\prime / \KK)$-invariant function on $\Xscr_{\KK^\prime}$, 
it descends to a function $g$ on $\Xscr_\KK$.
We have $\zero (g^\prime) = \nu^* (\zero (g))$,
and thus $\widetilde{\rho}_*^\prime ( \nu^* (\zero (g)) ) = d\, \underline{n} \in \Div (G)$.
By Claim~\ref{claim:surj-prins-2},
we obtain $\widetilde{\rho}_*(\zero (g)) = d\, \underline{n}$.
In conclusion, $\widetilde{L} := \zero (g)$ is a principal divisor on 
$\Xscr_\KK$ with $\widetilde{\rho}_* (\widetilde{L}) = d\, \underline{n}$.

Let $e > 2$ be another integer with $(e , d) = 1$.
When 
$p > 0$, we require that $(e, p) = 1$. By the above argument 
with $e$ in place of $d$, 
there exists a principal divisor
$\widetilde{M} \in \mathrm{Prin} (\Xscr_\KK)$ 
with $\widetilde{\rho}_* (\widetilde{M}) = e\, \underline{n}$.
We take integers $\alpha$ and $\beta$ such that  $\alpha d + \beta e = 1$,
and set $\widetilde{N} := \alpha \widetilde{L} + \beta \widetilde{M}$.
Then $\widetilde{N} \in \mathrm{Prin} (\Xscr_\KK)$ and 
$\widetilde{\rho}_*(\widetilde{L}) = \underline{n}$. 
This shows the theorem.  
\QED

%%%%%%%%%
% Bibliography %
%%%%%%%%%

\end{document}